\documentclass[12pt,reqno]{amsart}
  
 \usepackage{comment}
\usepackage{tikz}
\usepackage{mathrsfs}
\usepackage{latexsym}
\usepackage{graphicx}
\usepackage{amssymb}
 \graphicspath{{Figures/}}
 \usepackage[show]{ed}

\renewcommand{\epsilon}{\varepsilon}

\newcommand{\supp}{{\operatorname{supp\,}}}
\renewcommand{\phi}{\varphi}

\newtheorem{theo}{{\sc Theorem}}

\newtheorem{cor}{{\sc Corollary}}[section]

\newtheorem{prop}[cor]{{\sc Proposition}}
\numberwithin{equation}{section}

\newenvironment{rem}{\medskip\noindent{\it Remark:\/} }{\medskip}

\newtheorem{defn}[theo]{{\sc Definition}}

\usepackage{hyperref}

\title[Average estimate for Eigenfunctions]{Average estimate for Eigenfunctions along geodesics in the quantum completely integrable case}

\address{School of Mathematical Sciences\\ University of Electronic Science and Technology of China\\Chengdu 611731\\China}
\email{wawnwg123@163.com}
\author{Weiwei Wang And Xianchao Wu}
\address{School of Mathematics and Statistics, Wuhan University of Technology, Wuhan, Hubei, China}
\email{xianchao.wu@whut.edu.cn} 
\date{}

\begin{document}

\maketitle
\begin{abstract}
This paper investigates the upper bound of the integral of $L^2$-normalized joint eigenfunctions over geodesics in a two-dimensional quantum completely integrable system. For admissible geodesics, we rigorously establish an asymptotic decay rate of  $O(h^{1/2})$. This represents a polynomial improvement over the previously well known $O(1)$ bound.

\end{abstract}

\section{Introduction}
Let $(M,g)$ be a smooth compact 2-dimensional Riemannian manifold without boundary and $\{u_h\}$ be $L^2$-normalized eigenfunctions solving
\begin{equation}\label{laplace pde}
-h^2\Delta_g u_h=  u_h(x), \quad x\in M.
\end{equation}
The goal of this paper is to study the average oscillatory behavior of $u_h$ when restricted to a geodesic.

On a compact hyperbolic surface, using the Kuznecov formula, Good \cite{Goo} and Hejhal \cite{Hej} showed that if $\gamma$ is a periodic geodesic,  then
\begin{equation}\label{known}
\Big| \int_\gamma u_h(s)d s \Big|  \leq C_\gamma,
\end{equation}
where $ds$ denotes the arc length measure on $\gamma$. Zelditch \cite{Zel} generalized this result to the eigenfunctions restricting on a hypersurface inside a compact manifold, and Jung-Zelditch \cite{JZ} obtained an explicit polynomial bound  $O(h^{1/2}|\ln h|^{1/2})$ for a density-one subsequence. However, the latter estimate is not satisfied for all eigenfunctions. As shown by Eswarathasan \cite{Esw}, one expects 
\begin{equation}\label{expect}
\Big| \int_\gamma u_h(s)d s \Big|\sim h^{\frac{1}2}.
\end{equation}

The problem of improving the upper bound \eqref{known} for hyperbolic surfaces was raised and discussed in Pitt \cite{Pit} and Reznikov \cite{Rez}. As an analogy with the Lindel\"of conjecture for certain $L$-functions, Reznikov \cite{Rez} conjectured, for any given $\epsilon>0$, 
\begin{equation}\label{conjecture}
\Big| \int_\gamma u_h(s)d s \Big|\leq C_\epsilon h^{\frac{1}2 -\epsilon}
\end{equation} on a periodic geodesic $\gamma$. 
Working on surfaces of strictly negative curvature, Chen-Sogge \cite{CS} improved \eqref{known} to a $o(1)$ bound. Subsequently, Sogge-Xi-Zhang \cite{SXZ} obtained an $O(|\ln h|^{-1/2})$ bound under a relaxed curvature condition. Using Jacobi fields, Wyman \cite{Wym19, Wym20} generalized the results in \cite{CS} and \cite{SXZ} to curves satisfying certain curvature assumptions. The recent work of Canzani-Galkowski-Toth \cite{CGT} and Canzani-Galkowski \cite{CG} got $o(h^\frac{1-k}2)$ bounds for integrals over submanifolds of codimension $k$ under a weaker  geometric assumptions.

In this paper, we aim to consider the joint eigenfunctions of quantum completely integrable systems. Let $P_1(h): C^\infty(M)\to C^\infty(M)$ be a self-adjoint semiclassical pseudodifferential operator of order $m$ which is elliptic in the classical sense, i.e. $|p_1(x,\xi)|\geq c|\xi|^m-C$. We say that $P_1(h)$ is quantum completely integrable (QCI) if there exist functionally independent $h$-pseudodifferential operators $P_2(h), \dots, P_n(h)$ with property that
\begin{equation*}
[P_i(h), P_j(h)]=0; \quad i,j=1,\dots, n.
\end{equation*}
Given the joint eigenvalues $E(h)=(E_1(h), \dots, E_n(h))\in \mathbb{R}^n$ of $P_1(h), \dots, P_n(h)$ we denote an $L^2$-normalized joint eigenfunction with joint eigenvalue $E(h)$ by $u_{E,h}$. When the joint energy value $E$ is understood, we will sometimes abuse notation and simply write $u_h=u_{E,h}$.

The asymptotic bounds for the joint eigenfunctions $u_h$ of the QCI system are well studied. In \cite{Tot}, Toth  proved an asymptotic bounds for $\int_\gamma |u_h(s)|^2 ds$ on {\em generic} curves, and Galkowski-Toth \cite{GT2020} obtained a supremum bound for $\|u_h\|_{L^\infty}$ underunder a Morse type assumption. Recently,  Eswarathasan-Greenleaf-Keeler \cite{esw2024} established the pointwise Weyl law for the QCI system which satisfies the fiber rank $n$ condition. These results for the joint eigenfunctions $u_h$ of the QCI system are much better than the prediction in \cite{BGT} and the well-known H\"ormander bound $\|u\|_{L^\infty}=O(h^{(1-n)/2})$. 

In what follows, we attempt to obtain an asymptotic bounds for $\int_\gamma u_h(s)ds$ in the case where $\gamma$ is an {\em admissible} geodesic and $u_h$'s are {\em joint} eigenfunctions of the QCI system consisting of two commuting $h$-pseudodifferential operators $P_1(h)$ and $P_2(h)$ with joint eigenvalues $\big(E_1(h), E_2(h)\big)\in \text{Spec} P_1(h)\times \text{Spec} P_2(h)$. The admissible condition for the geodesic  is necessary.  Indeed, one can refer to the examples in Section \ref{example}, which show that  the bounds \eqref{known} cannot be improved if there is no assumption on the geodesic.

Before going on, let us explain the term {\em admissible}.

Let's assume that $p_1$ is of real principal type on the hypersurface $p_1^{-1}(E_1)$. This is equivalent to the following: if $E_1>0$ is a regular value of $p_1$, for any $(x,\xi)\in p_1^{-1}(E_1)$, the following inequality holds,
\begin{equation}\label{principal}
\partial_\xi p_1(x,\xi)\neq 0.
\end{equation}

Let $\pi: \, T^*M\to M$ be the canonical projection, we set
\begin{equation}
\mathcal{C}_\gamma:=\{(x,\xi)\in T^*M;\, p_1(x,\xi)=E_1,\, x\in \gamma\}=p_1^{-1}(E_1)\cap \pi^{-1}(\gamma).
\end{equation} 

\begin{rem}
In the homogeneous case where $p_1=|\xi|_g^2$ and $E_1=1$. One has that \eqref{principal} is automatically satisfied. 
\end{rem}

%\begin{defn}\label{def}
%Assume that a geodesic $\gamma$ is parametrized by its length $\tau$. We say that the $\mathbb{R}^2$-integrable system with moment map $\mathcal P=(p_1, p_2)$ is {\em admissible} along the geodesic $\gamma$ provided the condition \eqref{principal} is satisfied and following condition holds, 
%\begin{equation}\label{C1}
%\big|\partial_\tau p_2\,\vline_{\, C_\gamma\cap \{|p_2(x,\xi)-E_2|< \epsilon\}}\big|>C_\epsilon>0,
%\end{equation}
%for some sufficiently small positive constant $\epsilon$.
%\end{defn}
%
%{\color{red}
\begin{defn}\label{def}
Assume that a geodesic $\gamma$ is parametrized by arc length $\tau$. We say that the integrable system with moment map $\mathcal P=(p_1, p_2)$ is {\em admissible} along the geodesic $\gamma$ provided that the condition \eqref{principal} is satisfied and the following  holds:
\begin{equation}\label{C1}
\big|\partial_\tau p_2\,\vline_{\, \mathcal{C}_\gamma\cap \{|p_2(x,\xi)-E_2|< \epsilon\}}\big|>C_\epsilon>0,
\end{equation}
for some sufficiently small positive constant $\epsilon$.
\end{defn}

\begin{rem}
Geometrically condition \eqref{C1} tells that the energy surface $p_2(x,\xi)\,\vline_{\, \mathcal{C}_\gamma}=E$ intersects $T^*\gamma$ transversely if $|E-E_2|<\epsilon$.

\end{rem}

We shall call the geodesic $\gamma$ itself {\em admissible} when the conditions in Definition \ref{def} are satisfied. Now we can state the main theorem of this paper.

\begin{samepage}
\begin{theo}\label{main1}
Let $\{u_h\}$ be the $L^2$-normalized joint eigenfunctions of commuting operators $P_1(h)=-h^2\Delta_g+V$, where $V\in C^\infty(M,\mathbb{R})$, and $P_2(h)$ on a compact smooth Riemannian surface $(M,g)$ with $E_1(h)=E_1+O(h)$. Then for {\em admissible} geodesic $\gamma$ and $h\in (0, h_0]$ with a sufficiently small positive constant $h_0$, one has the following upper bound
\begin{equation}\label{main1 eq}
\Big| \int_\gamma u_h(s)d s \Big|  \leq C_\gamma\, h^{\frac{1}2}.
\end{equation}
where the constant $C_\gamma$ depends on the curve $\gamma$.
\end{theo}
\end{samepage}

\begin{rem}
The estimate \eqref{main1 eq} in Theorem \ref{main1} gives an explicit polynomial improvement over the well known bound \eqref{known} and is consistent with \eqref{conjecture}. Like the results in \cite{Tot, GT2020}, the above estimate is uniform over all energy values $\{E_2\in\mathbb{R};\, (E_1, E_2)\in \mathcal{P}(T^*M)\}$ if $\gamma$ is {\em admissilbe}. From the examples constructed in Section \ref{example}, one can see that the {\em admissibility} assumption of $\gamma$ is crucial. 

The methodology presented in current paper could be extended to general elliptic operators and smooth curves, although such an extension needs more intricate change of variables. We have chosen to state our results for the Schr\"odinger operator $P_1(h)$ and geodesic segment $\gamma$ in this work as the formulation retains its elegance and encapsulates the essential techniques.
\end{rem}

%\subsection{Outline of the paper}

%\noindent{\sc Acknowledgements:} 

\bigskip

\section{Proof of the Theorem \ref{main1}}
\subsection{Some reductions}
To prove \eqref{main1 eq}, we begin with a standard reduction. Without loss of generality, we always assume that the geodesic segment $\gamma$ has length one, and  the injectivity radius of $(M,g)$ is 10 or more and $\gamma$ is contained in one coordinate patch. We may work in local coordinates so that $\gamma$ is just
\begin{equation*}
\{(\tau,0): -1/2\leq \tau\leq 1/2\},
\end{equation*}
and we identify the point $\tau\in \gamma$ as its coordinate $(\tau,0)$.

Let us first fix a real-valued function $\chi\in \mathcal{S}(\mathbb{R})$ satisfying
\begin{equation}\label{spt}
\chi(0)=1\quad \text{and}\quad \hat\chi(t)=0,\quad |t|\geq\epsilon,
\end{equation}
here $\epsilon$ is a small positive constant.
Since $\chi\big(h^{-1}[P_1(h)-E_1(h)]\big)u_h=u_h$, in order to prove \eqref{main1 eq}, it suffices to show that
\begin{equation}\label{main eq}
\Big| \int_{-\frac{1}2}^{\frac{1}2} \big(\chi \big(h^{-1}[P_1(h)-E_1(h)]\big) f\big) (\tau)d \tau \Big|  \leq C_\gamma\, h^{\frac{1}2} \|f\|_{L^2(M)}.
\end{equation}

The joint spectrum of $P_1(h)$ (resp. $P_2(h)$) will be denoted by $\lambda_j^{(1)}(h)$ (resp. $\lambda_j^{(2)}(h)$) with $j=1,\,2,\,3,\dots$.
We note that the kernel of the operator $\chi \big(h^{-1}[P_1(h)-E_1(h)]\big)$ is given by
\begin{equation*}
\chi\big(h^{-1}[P_1(h)-E_1(h)]\big) (w,y)=\sum_j \chi(h^{-1}[\lambda_j^{(1)}(h)-E_1(h)]) u_j^h(w) \overline{u_j^h(y)},
\end{equation*}
here $u_j^h,\, j=1,2,3,\dots$ are the corresponding $L^2$-normalized joint eigenfunctions.

By the Cauchy-Schwarz's inequality, we would have \eqref{main eq} if we could obtain that
\begin{equation*}
\int_M \Big| \int_{-\frac{1}2}^{\frac{1}2} \sum_j \chi(h^{-1}[\lambda_j^{(1)}(h)-E_1(h)])u^h_j(\tau) \overline{u^h_j(y)} d\tau   \Big|^2 d y \leq C_\gamma h  .
\end{equation*}
By orthogonality of $\{u_j^h\}$, if we set $\rho(t)=(\chi(t))^2$, this is equivalent to showing that
\begin{equation}\label{upshot}
\Big|  \sum_j \rho(h^{-1}[\lambda_j^{(1)}(h)-E_1(h)])\int  \int u^h_j(\tau_1) \overline{u^h_j(\tau_2)} d\tau_1 d\tau_2 \Big| \leq C_\gamma h .
\end{equation}

%Now we fix another real-valued function $\rho_2\in \mathcal{S}(\mathbb{R})$ satisfying
%\begin{equation}\label{spt rho2}
%\rho_2(0)=1\quad \text{and}\quad \hat\rho_2(s)=0,\,\, s\geq C\epsilon.
%\end{equation}
We now claim that we would have \eqref{upshot} if we can prove the following:
\begin{prop}\label{propmain}
Under the condition of Theorem \ref{main1}, we have that, for $h\in (0, h_0]$,
\begin{align}\label{upshot5}
\sup_{\{E_2;\, (E_1,E_2)\in \mathcal{P}(T^*M)\}}\sum_j \rho(h^{-1}[\lambda_j^{(1)}(h)-E_1])\rho(h^{-1}[\lambda_j^{(2)}(h)-E_2]) \nonumber\\
\times \int_{-\frac{1}2}^{\frac{1}2}\int_{-\frac{1}2}^{\frac{1}2}  u^h_j (\tau_1) \overline{u^h_j(\tau_2)}  d\tau_1 d\tau_2 
\leq C_\gamma h.
\end{align}
where the constant $C_\gamma$ depends on the curve $\gamma$.
\end{prop}

Indeed, from \cite[Section 2]{Tot}, one knows that there exists a constant $C_2>0$ (independent of $j$ and $h$) such that for any $h\in(0, h_0]$, and $(\lambda_j^{(1)}(h), \lambda_j^{(2)}(h))\in \text{Spec}(P_1(h), P_2(h))$, with $|\lambda_j^{(1)}(h)-E_1|\leq C_1 h$,
\begin{equation*}
\inf_{E_2}|\lambda_j^{(2)}(h)-E_2|\leq C_2 h.
\end{equation*}
So, after taking $\epsilon$ in \eqref{spt} sufficiently small there exists a constant $C_3>0$ (independent of $j$ and $h$) such that for all $j\geq 1$ and $h\in(0, h_0]$,
\begin{equation*}
\sup_{\{E_2;\, (E_1,E_2)\in \mathcal{P}(T^*M),\, |\lambda_j^{(1)}(h)-E_1|\leq C_1 h\}}\rho(h^{-1}[\lambda_j^{(2)}(h)-E_2])\geq C_3>0.
\end{equation*}
Since the sum on the LHS of \eqref{upshot5} has non-negative terms, by restricting to $\{j;\,|\lambda_j^{(1)}(h)-E_1|\leq C_1 \}$ and (after taking $\epsilon$ small enough) using that $\rho(h^{-1}[\lambda_j^{(1)}(h)-E_1])\geq C_4>0$ for these eigenvalues, one gets that
\begin{equation*}
\sum_{\{j;\, |\lambda_j^{(1)}(h)-E_1|=O(h)\}} \int_{-\frac{1}2}^{\frac{1}2}  \int_{-\frac{1}2}^{\frac{1}2} u^h_j(\tau_1) \overline{u^h_j(\tau_2)} d\tau_1 d\tau_2 =O_\gamma(h),
\end{equation*}
which can deduce \eqref{upshot}.

\subsection{Integral representation of kernel}
Now we are going to prove \eqref{upshot5} in Proposition \ref{propmain}. Note that
\begin{align}
&\sum_j \rho(h^{-1}[\lambda_j^{(1)}(h)-E_1])\rho(h^{-1}[\lambda_j^{(2)}(h)-E_2])  u^h_j (x) \overline{u^h_j(y)} \nonumber \\
%=&\sum_j \rho(h^{-1}[\lambda_j^{(1)}(h)-E_1])\big(\rho(h^{-1}[P_2(h)-E_2])  u^h_j\big) (x) \overline{u^h_j(y)} \nonumber \\
%=& \rho(h^{-1}[P_2(h)-E_2])\sum_j \rho(h^{-1}[\lambda_j^{(1)}(h)-E_1])  u^h_j(x) \overline{u^h_j(y)} \nonumber \\
=& \Big(\rho \big(h^{-1}[P_2(h)-E_2]\big) \circ \rho \big(h^{-1}[P_1(h)-E_1]\big)\Big)(x,y).
\end{align}
%Here $\rho \big(h^{-1}[P_1(h)-E_1]\big)(\cdot,\cdot)$ is the kernel of the operator $\rho \big(h^{-1}[P_1(h)-E_1]\big)$.

\noindent Hence in order to prove \eqref{upshot5}, one needs to show that
\begin{equation}
\Big|\int\int \Big(\rho \big(h^{-1}[P_2(h)-E_2]\big)\circ \rho \big(h^{-1}[P_1(h)-E_1]\big) \Big)(\tau_1, \tau_2)  d\tau_1 d\tau_2\Big|\leq C_\gamma h.
\end{equation}

Note that
\begin{equation*}
\rho \big(h^{-1}[P_1(h)-E_1]\big)=\int \hat\rho(t) e^{\frac{i}h t(P_1(h)-E_1)}dt,
\end{equation*}
and we can write the Schwartz kernel of $e^{\frac{i}h t(P_1(h)-E_1)}$ in the form 
\begin{equation}\label{kerl1}
(2\pi h)^{-2} \int_{\mathbb{R}^2}e^{\frac{i}h [\varphi_1(t, w,\xi)-\left<y, \xi\right>-tE_1]} a_1(t, w,y,\xi;h) d\xi+ O(h^\infty),
\end{equation}
where $a_1\sim \sum_{j=0}^\infty a_{1,j} h^j$, $a_{1,j}\in C^\infty$, $a_{1,0}\geq C>0$ and $\varphi_1(t, w,\xi)$ solves the eikonal equation
\begin{equation}\label{ek1}
\partial_t \varphi_1(t, w,\xi)=p_1(w, \partial_w \phi_1), \quad \varphi_1(0, w, \xi)=\left<w,\xi \right>.
\end{equation}
Also
\begin{equation*}
\rho\big(h^{-1}[P_2(h)-E_2] \big)=\int_{\mathbb{R}}\hat{\rho}(s) e^{\frac{i}h s(P_2(h)-E_2)} ds.
\end{equation*}
Here $e^{\frac{i}h s(P_2(h)-E_2)}$ has a Schwartz kernel that is of the form
\begin{equation}\label{kerl2}
(2\pi h)^{-2} \int_{\mathbb{R}^2}e^{\frac{i}h [\varphi_2(s, x,\eta)-\left<w, \eta\right>-sE_2]} a_2(s, x,w,\eta;h) d\eta+ O(h^\infty),
\end{equation}
where $a_2\sim \sum_{j=0}^\infty a_{2,j} h^j$, $a_{2,j}\in C^\infty$, $a_{2,0}\geq C>0$ and $\varphi_2(s,x,\eta)$ solves the eikonal equation
\begin{equation}\label{ek2}
\partial_s \varphi_2(s, x, \eta)=p_2(x, \partial_x \phi_2), \quad \varphi_2(0, x, \eta)=\left<x,\eta \right>.
\end{equation}
From \eqref{ek1} and \eqref{ek2}, one can derive the following Taylor expansions for $\varphi_1(t, w,\xi)$ (resp. $\varphi_2(s, x,\eta)$) centered at $t=0$ (resp. $s=0$),
\begin{align*}
\varphi_1(t, w,\xi)&=w\xi+ tp_1(w,\xi)+O(t^2), \\
\varphi_2(s, x,\eta)&=x\eta+ sp_2(x,\eta)+O(s^2).
\end{align*}
Combining \eqref{kerl1} and \eqref{kerl2},  $\Big(\rho \big(h^{-1}[P_2(h)-E_2]\big) \circ \rho \big(h^{-1}[P_1(h)-E_1]\big)\Big)(x,y)$ is of the form
\begin{equation}\label{main kernel}
h^{- 4}\int \hat{\rho}(s) \hat{\rho}(t) \exp\left[\frac{i}h \Psi(s,t,x,w,y,\xi, \eta)\right]  
c(s,t,x,w,y,\xi, \eta;h) dw d\xi dt d\eta ds+O(h^\infty),
\end{equation}
where, $c\sim \sum_{j=0}^\infty c_{0,j} h^j$, $c_{0,j}\in C^\infty$, $c_{0,0}\geq C>0$ and
%\begin{align*}
%&\Phi(s,t,x,w,y,\xi, \eta) \nonumber \\
%=&\left<w-y,\xi \right>+\left<x-w, \eta\right>+t\big(p_1(w,\xi)-E_1\big)+s(p_2(x,\eta)-E_2)+t^2 f(w, \xi)  + s^2 g(x,\eta) \nonumber\\
%&  +O_{w,\xi}(t^3)+O_{x,\eta}(s^3).
%\end{align*}
\begin{align*}
&\Psi(s,t,x,w,y,\xi, \eta) \nonumber \\
=&\left<w-y,\xi \right>+\left<x-w, \eta\right>+t\big(p_1(w,\xi)-E_1\big)+s(p_2(x,\eta)-E_2)+O_{w,\xi}(t^2) +O_{x,\eta}(s^2).
\end{align*}
Up to an $O(h^\infty)$ term, we can write the amplitude in \eqref{main kernel} as
\begin{equation*}
\chi\big(p_1(w,\xi)-E_1\big) \chi\big(p_2(x,\eta)-E_2\big)c(s,t,x,w,y,\xi, \eta;h).
\end{equation*}
%Hence \eqref{main kernel1} is equal to
%\begin{equation}\label{main kernel}
%h^{- 4}\int \hat{\rho}(s) \hat{\rho}(t) \exp[\frac{i}h \Phi(s,t,x,w,y,\xi, \eta)]  
%c'(s,t,x,w,y,\xi, \eta;h) dw d\xi dt d\eta ds + O(h^\infty).
%\end{equation}

One can apply stationary phase in $(w, \xi)$-variables in \eqref{main kernel}. The critical point equations for $(w, \xi)$ are
\begin{align*}
\xi&=\eta+ t\partial_w p_1(w,\xi)+ O(t^2), \\
w&=y+t\partial_\xi p_1(w,\xi)+O(t^2).
\end{align*}
Since $\det(\Psi^{\prime\prime}_{w,\xi})=1+O(t)$, by stationary phase, \eqref{main kernel} equals
\begin{align*}
&h^{-2}\int \hat{\rho}(s) \hat{\rho}(t)  \exp\Big[\frac{i}h\Big(\left<x-y,\eta\right>+ t(p_1(y,\eta)-E_1\big)+s\big(p_2(x,\eta)-E_2\big)  \nonumber \\
&\qquad\qquad\qquad\qquad\quad+O_{x,\eta}(s^2)+O_{y,\eta}(t^2) \Big)\Big] c_1(s,t,x,y, \eta;h) d\eta dt ds+O(h^\infty) =: \tilde{I} (x,y; h),
\end{align*}
where $c_1\sim \sum_{j=0}^\infty c_{1,j} h^j$ and $c_{1,j}\in C^\infty$.

Substituting the points $\tau_1,\tau_2$ into $\tilde{I} (x,y; h)$, one only needs to prove that
\begin{equation}\label{upshot1}
\Big|\int \int   \tilde{I} (\tau_1, \tau_2; h) d\tau_1 d\tau_2\Big|\leq C_\gamma h.
\end{equation}

\subsection {\em Laplacian case.} If $P_1(h)=-h^2\Delta_g$ with $E_1=1$, %we use the geodesic polar coordinate about $\tau_2$ and 
 we can make the change of variables $\eta=(\eta_1,\eta_2)=r\Theta$, where $r=|\eta|$ and $\Theta=(\cos\theta, \sin\theta)\in S^{1}$ with $0\leq\theta<2\pi$. Then
\begin{align}\label{tilde I}
\tilde{I} (\tau_1, \tau_2; h)=& h^{-2}\int\int\int\int \hat{\rho}(s) \hat{\rho}(t)  \exp\left[\frac{i}h\Big(r \left<\tau_1-\tau_2, \Theta\right>+ t (r^2-1\big)+s\big(p_2(\tau_1, r\Theta)-E_2\big)\right. \nonumber\\ 
&\left. +O_{\tau_1,\eta}(t^2) +O_{\tau_2, \eta}(s^2)  \Big)\right] c_1(s,t, \tau_1, \tau_2, r\Theta;h) r dr d\Theta dt ds+O(h^\infty).
\end{align}
Recall that  the point $\tau_i\;(i=1,2)$ has coordinate $(\tau_i,0)\;(i=1,2)$. Then in \eqref{tilde I}
$$p_2(\tau_1, r\Theta)=p_2(\tau_1,0,\eta_1,\eta_2)=p_2(\tau_1,0,r\cos\theta,r\sin\theta)$$
and 
$$\langle\tau_1-\tau_2,\Theta\rangle=(\tau_1-\tau_2)\cos \theta.$$
To get \eqref{upshot1} from \eqref{tilde I}, we need the following three steps. \bigskip

{\it Step 1.} The critical point $\left(r_c, t_c\right)$ satisfies that
\begin{align}
\left<\tau_1-\tau_2, \Theta\right>+2rt  +s\Theta\cdot \partial_\eta p_2(\tau_1, r\Theta)+O_{\tau_1,\eta}(t^2) +O_{\tau_2, \eta}(s^2)=0, \label{r}\\
(r^2-1)+O_{\tau_1,\eta}(t)=0.  \nonumber
\end{align}
This gives that
\begin{align}
r_c&=r_c(\tau_1,\tau_2,s,\theta)=1+O(t_c), \label{r_c}\\
t_c&=t_c(\tau_1,\tau_2,s,\theta)=-\frac{1}{2}(\tau_1-\tau_2)\cos \theta +O(s).\label{E1}
\end{align}

%Noting that the determinant of Hessian of the phase function at $(r_c, t_c)$ is $1+O(s)$, we conclude that
%\cs\cs\cs
%\begin{align*}
%\Phi(s,\tau_1, \tau_2,\Theta) =(\tau_1-\tau_2)\cos \theta+s\left[p_2(\tau_1,0,\cos \theta,\sin\theta)-E_2\right]+O(s^2t_c)+O(t_c^2)+O(s^2).
%=&\left<\tau_1-\tau_2, \Theta\right>+s\big(p_2(\tau_1, \Theta)-E_2\big)+O(s t_c)+O(s^2)+O(t_c^2)
%\end{align*}
%\cs\cs\cs
Performing stationary phase at the critical point $\big(r_c, t_c\big)$ in \eqref{tilde I} gives,
$$\tilde{I} (\tau_1, \tau_2; h)=h^{-1}\int\int \hat{\rho}(s)\exp\left[\frac{i}h\Phi(s,\tau_1, \tau_2,\theta) \right] \cdot c_2(s, \tau_1, \tau_2, \theta; h) d\theta ds,$$
where $c_2\sim \sum_{j=0}^\infty c_{2,j} h^j$, $c_{2,j}\in C^\infty$ and
\begin{align*}
\Phi(s,\tau_1, \tau_2,\Theta)=&r_c\left<\tau_1-\tau_2, \Theta\right>+s\big(p_2(\tau_1, r_c \Theta)-E_2\big)+O(t_c^2)+O(s^2) \\
=&(\tau_1-\tau_2)\cos\theta+s\big(p_2(\tau_1, 0,\cos\theta,\sin\theta)-E_2\big)+O(s t_c)+O(s^2)+O(t_c^2)
\end{align*}
with the help of \eqref{r}, \eqref{r_c} and Taylor expansion in the last equality.\bigskip

{\it Step 2.} Take a new cut-off function $\chi\in C_0^\infty(\mathbb{R})$ with $\chi(\tau)=1$ if $\tau\in[0, 1]$ and $0\leq\chi(\tau)\leq1$. And set
\begin{equation}\label{chi1}
\chi_1(\theta)=\chi\left(\delta_0^{-1}\left(\theta-\frac{\pi}2\right) \right)+\chi\left(\delta_0^{-1}\left(\theta-\frac{3\pi}2\right) \right),
\end{equation}
where $\delta_0>0$ is a small constant. Insert this cut-off  into $ \tilde{I} (\tau_1, \tau_2; h)$:
\begin{align}
\tilde{I} (\tau_1, \tau_2; h)=&h^{-1}\int\int \hat{\rho}(s)\big(1-\chi_1(\theta)\big) \exp\left[\frac{i}h\Phi(s,\tau_1, \tau_2,\theta) \right] \cdot c_2(s, \tau_1, \tau_2, \theta; h) d\theta ds \nonumber \\
&+h^{-1}\int\int \hat{\rho}(s) \chi_1(\theta) \exp\left[\frac{i}h\Phi(s,\tau_1, \tau_2,\theta) \right] \cdot c_2(s, \tau_1, \tau_2, \theta; h) d\theta ds \nonumber \\
=&:\tilde{I}_1 (\tau_1, \tau_2; h)+\tilde{I}_2 (\tau_1, \tau_2; h).
\end{align}
For the term $\tilde{I}_1 (\tau_1,\tau_2; h)$,  notice that $|\cos\theta|>C>0$ if $\theta \in \supp\chi_1$. By taking the $\epsilon$ small enough in \eqref{spt}, one can use integration by parts in new variable $\tau:=\tau_1-\tau_2$ to get that
\begin{equation}
\left| \int\int \tilde{I}_1 (\tau_1,\tau_2; h) d\tau_1d\tau_2 \right|=O(h^\infty).
\end{equation}

Now it remains to to deal with the integral
$$\left|\int\int \tilde{I}_2 ( \tau_1 , \tau_2; h) d\tau_1 d\tau_2\right|.$$ 
One would like to apply stationary phase in $(\tau_1, \theta)$ to the integral
\begin{equation*}
 h^{-1}\int \hat{\rho}(s) \chi_1(\theta) \exp\left[\frac{i}h\Phi(s, \tau_1, \tau_2,\theta)\right] \cdot c_2(s, \tau_1, \tau_2, \theta; h) d\tau_1 d\theta.
\end{equation*}

Taking derivative with respect to $\tau_1$ in \eqref{E1} and choosing $\delta_0$ sufficiently small in \eqref{chi1}, one has that
\begin{equation}\label{t_c}
\frac{\partial{t_c}}{\partial{\tau_1}}=O(\epsilon).
\end{equation}

Notice that the critical point $\left(\tau_{1,c}(\tau_2, s), \theta_c(\tau_2, s) \right)$ is the solution of
\begin{align}
\cos\theta+ s\partial_{\tau_1}p_2+O(\epsilon)s+O(\epsilon)t_c=0, \label{tau_1}\\
-(\tau_1-\tau_2)\sin\theta+s\big(- \partial_{\eta_1} p_2 \cdot\sin\theta + \partial_{\eta_2} p_2\cdot \cos\theta \big)+O(\partial_\theta t_c)s+O(\partial_\theta t_c)t_c&=0, \label{theta}
\end{align}
with the help of \eqref{t_c} and reminding that $p_2=p_2(\tau_1,0, \eta_1, \eta_2)$. 

%\cs\cs Based on \eqref{derivoftau}, taking second order derivative with respect to $\tau_2$ in \eqref{theta}, one can get
%\begin{equation}\label{derivSecoftau}
%\partial_{\tau_2}^2\tau_{1,c}=O(s)\quad\mathrm{and}\quad\partial^2_{\tau_2}\theta_{c}=O(s),
%\end{equation}
%since $p_2$ is smooth when $\eta\in S^{1}$. 

Using \eqref{t_c},  the Hessian $\Phi^{\prime\prime}_{\tau_1, \theta}$ at the critical point $(\tau_{1,c}, \theta_c)$ is of the form
\begin{equation}\label{Hessian}
\begin{pmatrix}
O(\epsilon) & -\sin\theta_c+O(\epsilon) \\
-\sin\theta_c+O(\epsilon)  & O(1)
\end{pmatrix}.
\end{equation}
By taking the $\epsilon$ small enough in \eqref{spt} and $\delta_0$ sufficiently small \eqref{chi1}, one has that the above Hessian \eqref{Hessian} is non-degenerate. Then, performing the stationary phase in $(\tau_1, \theta)$ at the critical point $\left(\tau_{1,c}(\tau_2, s), \theta_c(\tau_2, s) \right)$ gives that
\begin{equation}\label{tau2sI2}
\int\int \tilde{I}_2 ( \tau_1 , \tau_2; h) d\tau_1 d\tau_2 =\int\int \exp\left[\frac{i}h  \tilde{\Phi}(\tau_2, s)\right]  \cdot c_3(\tau_2, s)d\tau_2 ds+O(h^\infty)
\end{equation}
with a new phase
\begin{align*}
\tilde{\Phi}(\tau_2, s) :=&(\tau_{1,c}-\tau_2)\cos\theta_c+s\left[p_2(\tau_{1,c}, 0, \cos\theta_c,\sin\theta_c)-E_2\right]+O(s t_c)+O(t_c^2)+O(s^2)
\end{align*}
 and $c_3\sim \sum_{j=0}^\infty c_{3,j} h^j$ with $c_{3,j}\in C^\infty$. Here $t_c$ take its value at $(\tau_{1,c},\tau_2,s,\theta_c)$.\bigskip

%We also notice that $\theta_0(\tau_1,\tau_2)\neq \frac{\pi}2\text{ or } \frac{3\pi}2$ if $\theta\in\supp{\chi_1}$, since $\partial_s\Phi(s, \tau_1, \tau_2,\Theta)=p_2(\tau_1,(\cos\theta,\sin\theta))-E_2+O(s)$ and $p_2(\tau_1, (\cos\theta,\sin\theta))\in \mathbb{R}\backslash[E_2-\epsilon_0, E_2+\epsilon_0]$ if $\theta\in\supp{\chi_1}$.

{\it Step 3.}  One needs to do integration by parts with respect to the variable $\tau_2$ to deal with the integration \eqref{tau2sI2}. 

First, some estimates are in order. Combining \eqref{E1} and  \eqref{tau_1} with taking the $\epsilon$ small enough in \eqref{spt}, one has that
\begin{equation}\label{ts}
t_c(\tau_{1,c}, \tau_2, s, \theta_c)=O(s).
\end{equation}
Hence 
\begin{equation}\label{t theta}
\frac{\partial{t_c}}{\partial \theta} (\tau_{1,c}, \tau_2, s, \theta_c)=O(s)
\end{equation}
and
\begin{equation*}
\tilde{\Phi}(\tau_2, s)=(\tau_{1,c}-\tau_2)\cos\theta_c+s\left[p_2(\tau_{1,c}, 0, \cos\theta_c,\sin\theta_c)-E_2\right]+O(s^2).
\end{equation*}

Furthermore, from \eqref{tau_1} and \eqref{theta}, one can get that
\begin{equation}\label{tauc thetac}
(\tau_{1,c}-\tau_2)\sin\theta_c=O(s) \quad\mathrm{and}\quad \cos\theta_c=O(s).
\end{equation}
Taking derivative with respect to $\tau_2$, one can conclude that
\begin{equation}\label{derivoftau}
(\partial_{\tau_2}\tau_{1,c}-1)\sin\theta=O(s)\quad\mathrm{and}\quad\partial_{\tau_2}\theta_{c}=O(s).
\end{equation}
Then the second and third order derivatives of $\tau_{1,c}$ and $\theta_{c}$ with respect to $\tau_2$ are
\begin{equation}\label{derivSecoftau}
\partial_{\tau_2}^k\tau_{1,c}=O(s)\quad\mathrm{and}\quad\partial^k_{\tau_2}\theta_{c}=O(s), \quad k=2,\, 3.
\end{equation}
%since $p_2$ is smooth when $\eta\in S^{1}$. 

Now using \eqref{theta}, \eqref{tauc thetac} and \eqref{derivoftau}, one can differentiate $ \tilde{\Phi}(\tau_2, s)$ with respect to $\tau_2$,
\begin{align}
& \partial_{\tau_2} \tilde{\Phi}(\tau_2, s) \nonumber\\
=&(\partial_{\tau_2}\tau_{1,c}-1)\cos\theta_c- (\tau_{1,c}-\tau_2)\sin\theta_c \partial_{\tau_2}\theta_c + s \partial_{\tau_2} p_2(\tau_{1,c},  \Theta_c)+O(s^2) \nonumber \\
=& O(s^2)+ s\Big( \partial_{\tau_2}\tau_{1,c}\cdot \partial_{\tau_1}p_2(\tau_{1,c},  \Theta_c) -  \partial_{\tau_2}\theta_c \cdot \big( \sin\theta_c \cdot \partial_{\eta_1} p_2(\tau_{1,c},  \Theta_c)+ \cos\theta_c \cdot \partial_{\eta_2} p_2(\tau_{1,c},  \Theta_c)\big)   \Big) \nonumber \\
=& s \partial_{\tau_2}\tau_{1,c}\cdot \partial_{\tau_1}p_2(\tau_{1,c},  \Theta_c) +O(s^2), \label{1d phi} 
\end{align}
with reminding that $p_2(\tau_{1,c},  \Theta_c)=p_2(\tau_{1,c},0,\cos  \theta_c,\sin  \theta_c)$.

Next, since the geodesic segment $\gamma$ is assumed to be {\em admissible}, from \eqref{C1},  it follows that
\begin{equation*}
\big| \partial_{\tau_1}p_2(\tau_{1,c},0,\cos  \theta_c,\sin  \theta_c) \big|>C>0.
\end{equation*}
So from \eqref{1d phi} and \eqref{derivoftau}, one can get the lower bound
\begin{equation}\label{positive}
\big|\partial_{\tau_2} \tilde{\Phi}(\tau_2, s)\big|>C|s|>0.
\end{equation}
Furthermore, by differentiating \eqref{1d phi} with respect to $\tau_2$ and using \eqref{derivSecoftau}, we have that
\begin{equation}\label{2d phi}
\big|\partial_{\tau_2}^k \tilde{\Phi}(\tau_2, s)\big|\leq C|s|, \quad k=2,\,3.
\end{equation}

Finally, one can separate the integration
\begin{align*}
\Big|\int \int   \tilde{I}_2 (\tau_2,s; h) d\tau_2 ds\Big|& \leq \Big|\int_{|s|\leq h} \int   \tilde{I}_2 (\tau_2,s; h) d\tau_2 ds\Big| + \Big|\int_{h<|s|<1} \int   \tilde{I}_2 (\tau_2,s; h) d\tau_2 ds\Big|.
\end{align*}

\noindent Obviously
\begin{equation}\label{last1}
\Big|\int_{|s|\leq h} \int   \tilde{I}_2 (\tau_2,s; h) d\tau_2 ds\Big|=O(h).
\end{equation}

\noindent For the second term, one can use integration by parts to get that
\begin{align}
\left|\int_{h<|s|<1} \int   \tilde{I}_2 (\tau_2,s; h) d\tau_2 ds\right|=& h \left|\int_{h<|s|<1} \int \frac{1}{\partial_{\tau_2} \tilde{\Phi}(\tau_2,s)}  \partial_{\tau_2}\left(\exp\left[\frac{i}h  \tilde{\Phi}(\tau_2, s)\right] \right) \cdot c_3(\tau_2, s) d\tau_2 ds\right| \nonumber\\
\leq& h \left|\int_{h<|s|<1} \int (\partial_{\tau_2} \tilde{\Phi})^{-1} \exp\left[\frac{i}h \tilde{ \Phi}(\tau_2, s)\right]  \cdot  \partial_{\tau_2}c_3(\tau_2, s) d\tau_2 ds\right| \nonumber\\
&+ h \left|\int_{h<|s|<1} \int (\partial_{\tau_2} \tilde{\Phi})^{-2}\, \partial^2_{\tau_2} \Phi \exp\left[\frac{i}h  \tilde{\Phi}(\tau_2, s)\right]  \cdot  c_3(\tau_2, s) d\tau_2 ds\right|. \nonumber
\end{align}

One more integration by parts, with the help of \eqref{positive} and \eqref{2d phi} we can get that
\begin{align}\label{last2}
&\left |\int_{h<|s|<1} \int   \tilde{I}_2 (\tau_2,s; h) d\tau_2 ds\right|  \nonumber \\
\leq& h^2 \left|\int_{h<|s|<1} \int (\partial_{\tau_2} \tilde{\Phi})^{-2} \partial_{\tau_2}\left(\exp\left[\frac{i}h \tilde{ \Phi}(\tau_2, s)\right] \right)  \cdot  \partial_{\tau_2}c_3(\tau_2, s) d\tau_2 ds\right| \nonumber\\
&+ h^2 \left|\int_{h<|s|<1} \int (\partial_{\tau_2} \tilde{\Phi})^{-3}\, \partial^2_{\tau_2} \Phi \,\partial_{\tau_2}\left(\exp\left[\frac{i}h  \tilde{\Phi}(\tau_2, s)\right]\right)  \cdot  c_3(\tau_2, s) d\tau_2 ds\right| \nonumber\\
\leq&  2 h^2 \left|\int_{h<|s|<1} \int (\partial_{\tau_2} \tilde{\Phi})^{-3} \partial^2_{\tau_2} \tilde{\Phi} \exp\left[\frac{i}h \tilde{ \Phi}(\tau_2, s)\right]   \cdot  \partial_{\tau_2}c_3(\tau_2, s) d\tau_2 ds\right| \nonumber\\
&+  h^2 \left|\int_{h<|s|<1} \int (\partial_{\tau_2} \tilde{\Phi})^{-2} \exp\left[\frac{i}h \tilde{ \Phi}(\tau_2, s)\right]   \cdot  \partial^2_{\tau_2}c_3(\tau_2, s) d\tau_2 ds\right| \nonumber\\
&+ h^2 \left|\int_{h<|s|<1} \int\left[ (\partial_{\tau_2} \tilde{\Phi})^{-3} \, \partial^3_{\tau_2} \Phi- 3 (\partial_{\tau_2} \tilde{\Phi})^{-4} \, \big(\partial^2_{\tau_2} \Phi\big)^2 \right] \,\exp\left[\frac{i}h  \tilde{\Phi}(\tau_2, s)\right]  \cdot  c_3(\tau_2, s) d\tau_2 ds\right| \nonumber\\
&+ h^2 \left|\int_{h<|s|<1} \int (\partial_{\tau_2} \tilde{\Phi})^{-3}\, \partial^2_{\tau_2} \Phi \exp\left[\frac{i}h  \tilde{\Phi}(\tau_2, s)\right]  \cdot  \partial_{\tau_2}c_3(\tau_2, s) d\tau_2 ds\right| \nonumber\\
\leq& C h^2 \left|\int_{h<|s|<1}\frac{1}{s^2} ds \right|=Ch.
\end{align}

 In conclusion, from \eqref{last1} and \eqref{last2} we finish the proof of \eqref{main1 eq} when $P_1(h)$ is Laplacian.

\subsection{\em Schr\"odinger case}
To treat the more general Schr\"odinger case, one can work in the Jacobi metric $g_E=(E-V)_+ g$ instead of Riemann metric $g$. Consequently, using geodesic normal coordinates in $g_E$ centered at $\tau_2$ like what we did in \eqref{tilde I}, one can follow the similar argument as in the homogeneous case.

\bigskip

\section{Examples}\label{example}
\subsection{Admissible geodesic on surface of revolution}

We consider the convex surface of revolution which has been well studied in \cite{Tot}. Using geodesic polar coordinates $(t, \varphi)\in [-1,1]\times [0, 2\pi]$, one can parametrize the convex surfaces of revolution, and so the metric $g$ is
$$g=dt^2+f^2(t)d\varphi^2,$$
where the profile function satisfies $f(-1)=f(1)=0$ and $f^{2}(t)$ is a non-negative Morse function with a single non-degenerate maximum at $t=t_0\in(-1,1)$. Then the principal symbols are
\begin{equation*}
p_1(t,\varphi; \xi_t, \xi_\varphi)=\xi_t^2+f^{-2}(t)\xi_\varphi^2,
\end{equation*}
and
\begin{equation*}
p_2(t,\varphi; \xi_t, \xi_\varphi)=\xi_\varphi.
\end{equation*}

{\it Case 1.} Consider the geodesic segment $\gamma=\{(t_0, \varphi);\, 0<\alpha\leq \varphi \leq \beta<2\pi\}$. One has that
\begin{equation*}
\mathcal{C}_{\gamma}=\{(t_0,\varphi; \xi_t, \xi_\varphi);\, \xi_t^2+f^{-2}(t_0)\xi_\varphi^2=1\}.
\end{equation*}

If $E_2=0$, then on the set $\mathcal{C}_\gamma\cap\{(t_0,\varphi; \xi_t, \xi_\varphi);\, |\xi_\varphi|< \delta\}$ where $\delta>0$ is sufficiently small, one has that $\partial_\varphi p_2=0$. 
Hence, in this case, $\gamma$ is {\em not admissible}.

{\it Case 2.} Consider another geodesic segment $\gamma'=\{(t, 0);\, -1<a\leq t \leq b<1\}$. One has that
\begin{equation*}
\mathcal{C}_{\gamma'}=\{(t,0; \xi_t, \xi_\varphi);\, \xi_t^2+f^{-2}(t)\xi_\varphi^2=1\}.
\end{equation*}

If $E_2\neq 0$, then on the set $\mathcal{C}_{\gamma'}\cap\{(t,0; \xi_t, \xi_\varphi);\, |\xi_\varphi-E_2|< \delta\}$ where $\delta>0$ is sufficiently small, by implicit function theorem one can use $t$ and $\xi_t$ to parametrize $\mathcal{C}_{\gamma'}$. Now
\begin{equation*}
p_2^2(t,\xi_t)=f^2(t)(1-\xi_t^2).
\end{equation*}
One can check that
\begin{equation*}
\partial_t (p_2^2)=2f(t)f'(t)(1-\xi_t^2)\neq 0 \quad \text{if } t\in[a,b]\backslash\{t_0\}.
\end{equation*}
Hence $\gamma'=\{(t, 0);\, -1<a\leq t \leq b<1\}$ is {\em admissible} if $a>t_0$ or $b<t_0$.

Next we focus on the case of standard sphere $\mathbb{S}^2$. In this case, $f(t)=\sqrt{1-t^2}$ and $t_0=0$.
\subsection{The  zonal harmonics}

Now we use the longitudinal coordinate $\theta\in [0, \pi]$ (i.e. $t=\cos\theta$), and latitudinal coordinates $\phi\in [0, 2\pi)$, then $\mathbb{S}^2\owns x=(\sin\theta \cos \phi, \sin \theta \sin \varphi, \cos \theta)$. One takes the zonal harmonics
$$u_k=C_k P_k(\cos\theta)$$
which solves
$$-\Delta u_k= k(k+1) u_k,$$
here $C_k\sim k^{\frac{1}2}$ is the $L^2$ normalisation factor, and $P_k$ is the Legendre polynomial. The eigenfrequency $h^{-1}=\sqrt{k(k+1)}\sim k$. 

This is a QCI system, and the corresponding $h$-Laplacian is $P_1(h):=-h^2\Delta$ with eigenvalue $E_1(h)=1+O(h)$. The commuting operator is $P_2(h)=hD_\varphi$ with $E_2(h)=0$. Consider a geodesic segment
$$\gamma_1=\{(\cos\varphi, \sin\varphi, 0);\, \varphi\in[0, \frac{\pi}3] \}.$$
From the case $1$ in Section $3.1$, we know that $\gamma_1$ is not admissible.

Next we need to estimate $\Big| \int_{\gamma_1} u_k ds \Big |$. From \cite[Theorem 8.21.2]{Sze}, one has
$$P_k(\cos\theta)=\sqrt{\frac{2}{\pi k \sin \theta}}\cos\big[(k+\frac{1}2)\theta-\frac{\pi}4\big]+O(k^{-3/2}).$$
Hence
$$P_k(0)=\sqrt{\frac{1}{\pi k }}+O(k^{-3/2}).$$

For such $\{u_k\}$, one can easily get
$$\Big| \int_{\gamma_1} u_k ds \Big |=\Big| \int_{0}^{\frac{\pi}3} u_k d\varphi \Big |\sim 1.$$
%Hence from above two examples, one can see that conditions \eqref{C1} and \eqref{C2} both are needed.
This implies that  the improved bounds in Theorem \ref{main1} may not hold for non-admissible geodesics. 

\subsection{The  tesseral spherical harmonics} Then we consider the $L^2$ normalized tesseral spherical harmonics
\begin{equation*}
u_\lambda(x)=Y_l^k(\theta,\varphi)=C_l^k P_l^k(\cos\theta)e^{ik\varphi},
\end{equation*}
where $P_l^k$ is the associated Legendre polynomial and $C_l^k$ is the $L^2$ normalisation factor and $h^{-1}=\lambda=(l(l+1))^{1/2}$. We take $l=2k$. In this case, the corresponding $h$-Laplacian $P_1(h):=-h^2\Delta$ is QCI with commuting $P_2(h)=hD_\varphi$ with $E_2(h)=\frac{1}2$. 

From \cite{BDWZ} there are two turning points at $\theta_0, \theta_1$ (with $\theta_0< \theta_1$) where $\sin\theta_k\sim 1/2$ for $k=0,1$. The eigenfunction oscillates in the allowed region $\{\theta\in(\theta_0, \theta_1)\}$ and decays exponentially in the forbidden region outside $\{\theta\in[\theta_0, \theta_1]\}$. In the transition region, which has a width of $O(h^{2/3})$, the eigenfunction behaves like the Airy function \cite[Proposition 4.3]{BDWZ} and has a magnitude of $O(h^{-1/6})$. 

Now consider the geodesic segment 
$$\gamma_2=\{(\sin\theta,0, \cos\theta); \theta\in[\theta_0-\delta_0, \theta_0]\},$$
where $\delta_0>0$ is a small positive constant. This segment is a longitudinal arc from the forbidden region to the caustic point. From the case $2$ in Section $3.1$, it follows that $\gamma_2$ is {\em admissible}. One has that
\begin{equation*}
\Big| \int_{\gamma_2} u_\lambda(x) ds \Big|=O(h^{-1/6}\cdot h^{2/3})=O(h^{1/2}).
\end{equation*}
These bounds are consistent with the improved bounds in Theorem \ref{main1}.

\bigskip
\bibliography{reference}

\begin{thebibliography}{BDWZ12}

\bibitem[BDWZ12]{BDWZ}
Nicolas Burq, Semyon Dyatlov, Rachel Ward, and Maciej Zworski.
\newblock Weighted eigenfunction estimates with applications to compressed
  sensing.
\newblock {\em SIAM J. Math. Anal.}, 44(5):3481--3501, 2012.

\bibitem[BGT07]{BGT}
N.~Burq, P.~G\'{e}rard, and N.~Tzvetkov.
\newblock Restrictions of the {L}aplace-{B}eltrami eigenfunctions to
  submanifolds.
\newblock {\em Duke Math. J.}, 138(3):445--486, 2007.

\bibitem[CG19]{CG}
Yaiza Canzani and Jeffrey Galkowski.
\newblock On the growth of eigenfunction averages: microlocalization and
  geometry.
\newblock {\em Duke Math. J.}, 168(16):2991--3055, 2019.

\bibitem[CGT18]{CGT}
Yaiza Canzani, Jeffrey Galkowski, and John~A. Toth.
\newblock Averages of eigenfunctions over hypersurfaces.
\newblock {\em Comm. Math. Phys.}, 360(2):619--637, 2018.

\bibitem[CS15]{CS}
Xuehua Chen and Christopher~D. Sogge.
\newblock On integrals of eigenfunctions over geodesics.
\newblock {\em Proc. Amer. Math. Soc.}, 143(1):151--161, 2015.

\bibitem[ESB24]{esw2024}
Greenleaf~A. Eswarathasan~S. and Keeler B.
\newblock Pointwise weyl laws for quantum completely integrable systems.
\newblock {\em arXiv:2411.10401 [math.AP]}, 2024.

\bibitem[Esw16]{Esw}
Suresh Eswarathasan.
\newblock Expected values of eigenfunction periods.
\newblock {\em J. Geom. Anal.}, 26(1):360--377, 2016.

\bibitem[Goo83]{Goo}
Anton Good.
\newblock {\em Local analysis of {S}elberg's trace formula}, volume 1040 of
  {\em Lecture Notes in Mathematics}.
\newblock Springer-Verlag, Berlin, 1983.

\bibitem[GT20]{GT2020}
Jeffrey Galkowski and John~A Toth.
\newblock Pointwise bounds for joint eigenfunctions of quantum completely
  integrable systems.
\newblock {\em Communications in Mathematical Physics}, 375(2):915--947, 2020.

\bibitem[Hej82]{Hej}
Dennis~A. Hejhal.
\newblock Sur certaines s\'{e}ries de {D}irichlet associ\'{e}es aux
  g\'{e}od\'{e}siques ferm\'{e}es d'une surface de {R}iemann compacte.
\newblock {\em C. R. Acad. Sci. Paris S\'{e}r. I Math.}, 294(8):273--276, 1982.

\bibitem[JZ16]{JZ}
Junehyuk Jung and Steve Zelditch.
\newblock Number of nodal domains and singular points of eigenfunctions of
  negatively curved surfaces with an isometric involution.
\newblock {\em J. Differential Geom.}, 102(1):37--66, 2016.

\bibitem[Pit08]{Pit}
Nigel J.~E. Pitt.
\newblock A sum formula for a pair of closed geodesics on a hyperbolic surface.
\newblock {\em Duke Math. J.}, 143(3):407--435, 2008.

\bibitem[Rez15]{Rez}
Andre Reznikov.
\newblock A uniform bound for geodesic periods of eigenfunctions on hyperbolic
  surfaces.
\newblock {\em Forum Math.}, 27(3):1569--1590, 2015.

\bibitem[SXZ17]{SXZ}
Christopher~D. Sogge, Yakun Xi, and Cheng Zhang.
\newblock Geodesic period integrals of eigenfunctions on {R}iemannian surfaces
  and the {G}auss-{B}onnet theorem.
\newblock {\em Camb. J. Math.}, 5(1):123--151, 2017.

\bibitem[Sze75]{Sze}
G\'{a}bor Szeg\H{o}.
\newblock {\em Orthogonal polynomials}.
\newblock American Mathematical Society Colloquium Publications, Vol. XXIII.
  American Mathematical Society, Providence, RI, fourth edition, 1975.

\bibitem[Tot09]{Tot}
John~A. Toth.
\newblock {$L^2$}-restriction bounds for eigenfunctions along curves in the
  quantum completely integrable case.
\newblock {\em Comm. Math. Phys.}, 288(1):379--401, 2009.

\bibitem[Wym19]{Wym19}
Emmett~L. Wyman.
\newblock Looping directions and integrals of eigenfunctions over submanifolds.
\newblock {\em J. Geom. Anal.}, 29(2):1302--1319, 2019.

\bibitem[Wym20]{Wym20}
Emmett~L. Wyman.
\newblock Explicit bounds on integrals of eigenfunctions over curves in
  surfaces of nonpositive curvature.
\newblock {\em J. Geom. Anal.}, 30(3):3204--3232, 2020.

\bibitem[Zel92]{Zel}
Steven Zelditch.
\newblock Kuznecov sum formulae and {S}zeg\H{o} limit formulae on manifolds.
\newblock {\em Comm. Partial Differential Equations}, 17(1-2):221--260, 1992.

\end{thebibliography}
\bibliographystyle{alpha}

\end{document}